\documentclass{article}
\usepackage{amssymb,amsfonts,amsmath,amsthm}
\usepackage{hyperref}
\usepackage{parskip}
\usepackage{setspace}
\usepackage{fancyhdr}
\usepackage[dvipsnames]{xcolor}
\usepackage[numbers]{natbib}
\usepackage[bottom]{footmisc}
\usepackage[utf8]{inputenc}
\usepackage[english]{babel}
\usepackage{xcolor}
\usepackage{geometry}
\usepackage{graphicx}
\newcommand\myshade{85}
\colorlet{mylinkcolor}{blue}
\colorlet{mycitecolor}{red}
\colorlet{myurlcolor}{Aquamarine}
\hypersetup{
	linkcolor  = mylinkcolor,
	citecolor  = mycitecolor,
	urlcolor   = myurlcolor!\myshade!black,
	colorlinks = true,
}
\title{Tur\'an-type inequaities for generalized polygamma function}

\author{Omprakash Atale${}^{1}$\footnote{${}^{1}$ Department of Mathematics, 
Khandesh College Education Society's. Moolji Jaitha College."An Autonomous College Affiliated to K.B.C. North Maharashtra University, Jalgaon."Maharashtra, India.}} 
\date{July 1, 2022}

\begin{document}

\maketitle

\begin{abstract}
Inspired by the work of C. Mortici [1] and A. Laforgia et. al [2] we have established some new Turán-type inequalities for k-polygamma function and p-k-polygamma function. 
\end{abstract}

\section{Introduction}
In 1950, a mathematician named Pal Turan derived the following inequality for Legendre polynomial [1]:
\begin{equation} \label{GrindEQ__1_1_} 
P_{n} (x)^{2} >P_{n-1} \left(x\right)P_{n+1} \left(x\right) 
\end{equation} 
valid for -1 $\mathrm{<}$ x $\mathrm{<}$ 1. This inequality is now known as Turan's inequality and have wide applications in wide areas of mathematics such as complex analysis, number theory, combinatorics, theory of mean-values or statistics and control theory. Turan's inequalities for Hermite polynomial and Chebyshev polynomials (-1 $\mathrm{<}$ x $\mathrm{<}$ 1) are
\begin{equation} \label{GrindEQ__1_2_} 
H_{n} \left(x\right)^{2} -H_{n-1} \left(x\right)H_{n+1} \left(x\right)=\left(n-1\right)!\sum _{i=0}^{n-1}\frac{2^{n-i} }{i!}  H_{i} \left(x\right)^{2} >0 
\end{equation} 
and
\begin{equation} \label{GrindEQ__1_3_} 
T_{n} \left(x\right)^{2} -T_{n-1} \left(x\right)T_{n+1} \left(x\right)=1-x^{2} >0 
\end{equation} 
respectively. Recently, C. Mortici [2] proved some new type of Turan's inequalities for polygama function using some new results derived previously by Laforgia and Natalini [3]. In this paper, we are going to do something similar for the generalizations of polygama function, namely, the \textit{k-}polygamma function and the \textit{p-k-}polygamma function. First we introduce what \textit{k-}polygamma function and \textit{p-k-}polygamma function are and then we proceed towards our main results. Diaz and Pariguan [4] introduced the following generalization of gamma function known as the \textit{k-}gamma function:
\begin{equation} \label{GrindEQ__1_4_} 
\Gamma _{k} \left(x\right)=\int _{0}^{\infty }t^{x-1} e^{-\frac{t^{k} }{k} }  dt 
\end{equation} 
valid for Re(x) $\mathrm{>}$ 0 and k $\mathrm{>}$ 0. The corresponding polygamma function, which we call the \textit{k-}polygamma function is given by the definition 
\begin{equation} \label{GrindEQ__1_5_} 
\psi _{k}^{\left(m\right)} \left(x\right)=\left(-1\right)^{m+1} m!\sum _{n=0}^{\infty }\frac{1}{\left(nk+x\right)^{m+1} }  =\left(-1\right)^{m+1} \int _{0}^{\infty }\frac{t^{m} e^{-xt} }{1-e^{-kt} }  dt.                                     
\end{equation} 
Soon after, Gehlot [5] introduced the following generalization of the \textit{k-}gamma function:
\begin{equation} \label{GrindEQ__1_6_} 
{}_{p} \Gamma _{k} \left(x\right)=\int _{0}^{\infty }t^{x-1} e^{-\frac{t^{k} }{p} }  dt 
\end{equation} 
valid for Re(x) $\mathrm{>}$ 0 and p, k $\mathrm{>}$ 0. The ordinary gamma function, \textit{k-}gamma function and the \textit{p-k-}gamma function are related to each-other as follows:
\begin{equation} \label{GrindEQ__1_7_} 
{}_{p} \Gamma _{k} \left(x\right)=\left(\frac{p}{k} \right)^{\frac{x}{k} } \Gamma _{k} \left(x\right)=\frac{p^{\frac{x}{k} } }{k} \Gamma \left(\frac{x}{k} \right).            
\end{equation} 
\section{Some results using Holder's inequaities}
  Let $\xi $ and $\tilde{\xi }$ be two non-negative integrable functions over the range of $[a,b]$. Let $p,q>0$such that $p^{-1} +q^{-1} =1$. Then, we define Holder's inequality as
\begin{equation} \label{GrindEQ__2_1_} 
\left(\int _{a}^{b}\xi ^{p} \left(t\right)dt \right)^{\frac{1}{p} } \left(\int _{a}^{b}\tilde{\xi }^{q} \left(t\right)dt \right)^{\frac{1}{q} } \ge \int _{a}^{b}\xi \left(t\right)\tilde{\xi }\left(t\right)dt  
\end{equation} 
For the particular case of $p=q=2$, the above inequality is reduced to the Cauchy-Schwarz inequality. Let $\xi =g^{\frac{1}{p} } \left(t\right)f^{\frac{m}{p} } \left(t\right)$ and $\tilde{\xi }=g^{\frac{1}{q} } \left(t\right)f^{\frac{n}{q} } \left(t\right)$to get
\begin{equation} \label{GrindEQ__2_2_} 
\left(\int _{a}^{b}g\left(t\right)f^{m} \left(t\right)dt \right)^{\frac{1}{p} } \left(\int _{a}^{b}g\left(t\right)f^{n} \left(t\right)dt \right)^{\frac{1}{q} } \ge \int _{a}^{b}g\left(t\right)f^{\frac{m}{p} +\frac{n}{q} } \left(t\right)dt .                                  
\end{equation} 
\textbf{Theorem 1.} \textit{For every $p,q>0$ with $p^{-1} +q^{-1} =1$ and $m,n\ge 1$ such that $\frac{m}{p} +\frac{n}{q} $ is an integer, we have}
\begin{equation} \label{GrindEQ__2_3_} 
\left(\psi _{k}^{\left(m\right)} \left(x\right)\right)^{\frac{1}{p} } .\left(\psi _{k}^{\left(n\right)} \left(x\right)\right)^{\frac{1}{q} } \ge \psi _{k}^{\left(\frac{m}{p} +\frac{n}{q} \right)} \left(x\right). 
\end{equation} 
\textit{Proof:} Let $g\left(t\right)=\frac{e^{-xt} }{1-e^{-kt} } $, $f\left(t\right)=t$and $a=0,b=+\infty $in Eqn. \eqref{GrindEQ__2_2_} to get
\begin{equation} \label{GrindEQ__2_4_} 
\left(\int _{0}^{\infty }\frac{t^{m} e^{-xt} }{1-e^{-kt} } dt \right)^{\frac{1}{p} } \left(\int _{0}^{\infty }\frac{t^{n} e^{-xt} }{1-e^{-kt} } dt \right)^{\frac{1}{q} } \ge \int _{0}^{\infty }\frac{t^{\frac{m}{p} +\frac{n}{q} } e^{-xt} }{1-e^{-kt} } dt .                                               
\end{equation} 
Using the integral representation of $\psi _{k}^{\left(m\right)} \left(x\right)$ from Eqn. \eqref{GrindEQ__1_5_}, the desired result readily follows. 

          In [4], we have the following proposition for $a\in {\rm R}$
\begin{equation} \label{GrindEQ__2_5_} 
\Gamma _{k} \left(x\right)=a^{\frac{x}{k} } \int _{0}^{\infty }t^{x-1} e^{-\frac{t^{k} }{k} a}  dt.                                                                  
\end{equation} 
From this, we can get the following integral representation of the \textit{k-}zeta function
\begin{equation} \label{GrindEQ__2_6_} 
\zeta _{k} \left(x\right)\equiv \zeta \left(\frac{x}{k} \right)=\frac{1}{\Gamma _{k} \left(x\right)} \int _{0}^{\infty }\frac{t^{x-1} }{e^{\frac{t^{k} }{k} } -1}  dt. 
\end{equation} 
\textbf{Theorem 2.} \textit{For every $p,q>0$ with $p^{-1} +q^{-1} =1$ and $m,n\ge 1$ such that $\frac{m}{p} +\frac{n}{q} $ is an integer, we have}
\begin{equation} \label{GrindEQ__2_7_} 
\left(\zeta _{k} \left(m+1\right)\right)^{\frac{1}{p} } \cdot \left(\zeta _{k} \left(n+1\right)\right)^{\frac{1}{q} } \ge \frac{\Gamma _{k} \left(\frac{m}{p} +\frac{n}{q} +1\right)}{\Gamma _{k} \left(m+1\right)^{\frac{1}{p} } \Gamma _{k} \left(n+1\right)^{\frac{1}{q} } } \zeta _{k} \left(\frac{m}{p} +\frac{n}{q} +1\right).                          
\end{equation} 
\textit{Proof:} Let replace $x$ with $x+1$in Eqn. \eqref{GrindEQ__2_6_} to get
\begin{equation} \label{GrindEQ__2_8_} 
\zeta _{k} \left(x+1\right)\equiv \zeta \left(\frac{x+1}{k} \right)=\frac{1}{\Gamma _{k} \left(x+1\right)} \int _{0}^{\infty }\frac{t^{x} }{e^{\frac{t^{k} }{k} } -1}  dt. 
\end{equation} 
Now, let$g\left(t\right)=\frac{1}{e^{\frac{t^{k} }{k} } -1} $, $f\left(t\right)=t$ and $a=0,b=+\infty $ in Eqn. \eqref{GrindEQ__2_2_} to get
\begin{equation} \label{GrindEQ__2_9_} 
\left(\int _{0}^{\infty }\frac{t^{m} }{e^{\frac{t^{k} }{k} } -1} dt \right)^{\frac{1}{p} } \left(\int _{0}^{\infty }\frac{t^{n} }{e^{\frac{t^{k} }{k} } -1} dt \right)^{\frac{1}{q} } \ge \int _{0}^{\infty }\frac{t^{\frac{m}{p} +\frac{n}{q} } }{e^{\frac{t^{k} }{k} } -1} dt ,                                                  
\end{equation} 
Using Eqn. \eqref{GrindEQ__2_6_}, we get
\begin{equation}
\left(\zeta _{k} \left(m+1\right)\Gamma _{k} \left(m+1\right)\right)^{\frac{1}{p} } \cdot \left(\zeta _{k} \left(n+1\right)\Gamma _{k} \left(n+1\right)\right)^{\frac{1}{q} } \ge \zeta _{k} \left(\frac{m}{p} +\frac{n}{q} +1\right)\Gamma _{k} \left(\frac{m}{p} +\frac{n}{q} +1\right)    
\end{equation}

  This completes out proof.

            In a similar manner, we get the following theorem for \textit{p-k-}gamma function. Consider the following proposition from [5]:
\begin{equation} \label{GrindEQ__2_11_} 
{}_{p} \Gamma _{k} \left(x\right)=a^{\frac{x}{k} } \int _{0}^{\infty }t^{x-1} e^{-\frac{t^{k} }{p} a}  dt.                                                                 
\end{equation} 
From this, we can get the following integral representation of the \textit{k-}zeta function
\begin{equation} \label{GrindEQ__2_12_} 
{}_{p} \zeta _{k} \left(x\right)=\frac{1}{{}_{p} \Gamma _{k} \left(x\right)} \int _{0}^{\infty }\frac{t^{x-1} }{e^{\frac{t^{k} }{p} } -1}  dt. 
\end{equation}

  \textbf{Theorem 3.} \textit{For every $r,s>0$ with $r^{-1} +s^{-1} =1$ and $m,n\ge 1$ such that $\frac{m}{r} +\frac{n}{s} $ is an integer, we have}
\begin{equation} \label{GrindEQ__2_13_} 
\left({}_{p} \zeta _{k} \left(m+1\right)\right)^{\frac{1}{r} } \cdot \left({}_{p} \zeta _{k} \left(n+1\right)\right)^{\frac{1}{s} } \ge \frac{{}_{p} \Gamma _{k} \left(\frac{m}{r} +\frac{n}{s} +1\right)}{{}_{p} \Gamma _{k} \left(m+1\right)^{\frac{1}{r} } {}_{p} \Gamma _{k} \left(n+1\right)^{\frac{1}{s} } } \cdot _{p} \zeta _{k} \left(\frac{m}{r} +\frac{n}{s} +1\right).                         
\end{equation}

  \textit{Proof: }Let replace $x$with$x+1$in Eqn. \eqref{GrindEQ__2_12_} to get
\begin{equation} \label{GrindEQ__2_14_} 
{}_{p} \zeta _{k} \left(x+1\right)=\frac{1}{{}_{p} \Gamma _{k} \left(x+1\right)} \int _{0}^{\infty }\frac{t^{x} }{e^{\frac{t^{k} }{p} } -1}  dt. 
\end{equation} 
Now, let$g\left(t\right)=\frac{1}{e^{\frac{t^{k} }{p} } -1} $, $f\left(t\right)=t$and $a=0,b=+\infty $ in Eqn. \eqref{GrindEQ__2_2_} to get
\begin{equation} \label{GrindEQ__2_15_} 
\left(\int _{0}^{\infty }\frac{t^{m} }{e^{\frac{t^{k} }{p} } -1} dt \right)^{\frac{1}{r} } \left(\int _{0}^{\infty }\frac{t^{n} }{e^{\frac{t^{k} }{p} } -1} dt \right)^{\frac{1}{s} } \ge \int _{0}^{\infty }\frac{t^{\frac{m}{r} +\frac{n}{s} } }{e^{\frac{t^{k} }{p} } -1} dt ,                                                 
\end{equation} 
Using Eqn. \eqref{GrindEQ__2_12_}, we get
\begin{equation}
\left({}_{p} \zeta _{k} \left(m+1\right)\cdot _{p} \Gamma _{k} \left(m+1\right)\right)^{\frac{1}{r} } \cdot \left({}_{p} \zeta _{k} \left(n+1\right)\cdot _{p} \Gamma _{k} \left(n+1\right)\right)^{\frac{1}{s} } \ge _{p} \zeta _{k} \left(\frac{m}{r} +\frac{n}{s} +1\right)\cdot _{p} \Gamma _{k} \left(\frac{m}{r} +\frac{n}{s} +1\right)    
\end{equation}

  This completes out proof.
\section{Turan-type inequalities}
  \textbf{Theorem 4.} \textit{We have}

  \textit{i)} 
\begin{equation} \label{GrindEQ__3_1_} 
0\le \Gamma _{k}^{\left(n-1\right)} \left(x\right)\Gamma _{k}^{\left(n+1\right)} \left(x\right)-\left(\Gamma _{k}^{\left(n\right)} \left(x\right)\right)^{2}  
\end{equation} 
\textit{and}

  \textit{ii)}
\begin{equation} \label{GrindEQ__3_2_} 
0\le _{p} \Gamma _{k}^{\left(n-1\right)} \left(x\right)\cdot _{p} \Gamma _{k}^{\left(n+1\right)} \left(x\right)-\left({}_{p} \Gamma _{k}^{\left(n\right)} \left(x\right)\right)^{2}  
\end{equation} 
\textit{Proof: }Take the n${}^{th}$ derivative of \textit{k-}gamma function to get
\begin{equation} \label{GrindEQ__3_3_} 
\Gamma _{k}^{\left(n\right)} \left(x\right)=\int _{0}^{\infty }t^{x-1} e^{-\frac{t^{k} }{k} } \log ^{n} t dt 
\end{equation} 
Now, using Cauchy-Schwarz inequality we get
\begin{equation} \label{GrindEQ__3_4_} 
\left(\Gamma _{k}^{\left(n\right)} \left(x\right)\right)^{2} =\left(\int _{0}^{\infty }\left[e^{-\frac{t^{k} }{k} } t^{x-1} \log \left(t\right)^{n-1} \right]^{\frac{1}{2} } \left[e^{-\frac{t^{k} }{k} } t^{x-1} \log \left(t\right)^{n+1} \right]^{\frac{1}{2} }  dt\right)^{2}  
\end{equation} 
\begin{equation} \label{GrindEQ__3_5_} 
\le \int _{0}^{\infty }e^{-\frac{t^{k} }{k} } t^{x-1} \log \left(t\right)^{n-1}  dt\int _{0}^{\infty }e^{-\frac{t^{k} }{k} } t^{x-1} \log \left(t\right)^{n+1}  dt=\Gamma _{k}^{\left(n-1\right)} \left(x\right)\Gamma _{k}^{\left(n+1\right)} \left(x\right) 
\end{equation} 
which implies
\begin{equation} \label{GrindEQ__3_6_} 
0\le \Gamma _{k}^{\left(n-1\right)} \left(x\right)\Gamma _{k}^{\left(n+1\right)} \left(x\right)-\left(\Gamma _{k}^{\left(n\right)} \left(x\right)\right)^{2} .                                                       
\end{equation} 
Similarly, for \textit{p-k-}gamma function we have
\begin{equation} \label{GrindEQ__3_7_} 
0\le _{p} \Gamma _{k}^{\left(n-1\right)} \left(x\right)\cdot _{p} \Gamma _{k}^{\left(n+1\right)} \left(x\right)-\left({}_{p} \Gamma _{k}^{\left(n\right)} \left(x\right)\right)^{2} .                                                  
\end{equation} 
\textbf{}

  \textbf{Theorem 5.} \textit{For $x>0$ and even integers $n\ge l\ge 0$, we have}
\begin{equation} \label{GrindEQ__3_8_} 
e^{\Gamma _{k}^{\left(n-l\right)} \left(x\right)} \cdot e^{\Gamma _{k}^{\left(n+l\right)} \left(x\right)} \ge \left(e^{\Gamma _{k}^{\left(n\right)} \left(x\right)} \right)^{2} .                                                                
\end{equation} 
\textit{Proof: }Using Eqn. \eqref{GrindEQ__3_3_}, estimate the expression
\begin{equation} \label{GrindEQ__3_9_} 
\frac{\Gamma _{k}^{\left(n-l\right)} \left(x\right)+\Gamma _{k}^{\left(n+l\right)} \left(x\right)}{2} -\Gamma _{k}^{\left(n\right)} \left(x\right) 
\end{equation} 
\begin{equation} \label{GrindEQ__3_10_} 
=\frac{1}{2} \left(\int _{0}^{\infty }t^{x-1} e^{-\frac{t^{k} }{k} } \log ^{n-l} t dt+\int _{0}^{\infty }t^{x-1} e^{-\frac{t^{k} }{k} } \log ^{n+l} t dt\right)-\int _{0}^{\infty }t^{x-1} e^{-\frac{t^{k} }{k} } \log ^{n} t dt 
\end{equation} 
\begin{equation} \label{GrindEQ__3_11_} 
=\frac{1}{2} \int _{0}^{\infty }\left(\frac{1}{\log ^{l} t} +\log ^{l} t-2\right)t^{x-1} e^{-\frac{t^{k} }{k} } \log ^{n} t dt\ge 0 
\end{equation} 
\begin{equation} \label{GrindEQ__3_12_} 
\frac{\Gamma _{k}^{\left(n-l\right)} \left(x\right)+\Gamma _{k}^{\left(n+l\right)} \left(x\right)}{2} \ge \Gamma _{k}^{\left(n\right)} \left(x\right) 
\end{equation} 
Exponentiating the above inequality yields the desired result.

  \textbf{Theorem 6.} \textit{For $x>0$ and even integers $n\ge l\ge 0$, we have}
\begin{equation} \label{GrindEQ__3_13_} 
e^{{}_{p} \Gamma _{k}^{\left(n-l\right)} \left(x\right)} \cdot e^{{}_{p} \Gamma _{k}^{\left(n+l\right)} \left(x\right)} \ge \left(e^{{}_{p} \Gamma _{k}^{\left(n\right)} \left(x\right)} \right)^{2} .                                                               
\end{equation} 
\textit{Proof: }Take the n${}^{th}$ derivative of \textit{k-}gamma function to get
\begin{equation} \label{GrindEQ__3_14_} 
{}_{p} \Gamma _{k}^{\left(n\right)} \left(x\right)=\int _{0}^{\infty }t^{x-1} e^{-\frac{t^{k} }{p} } \log ^{n} t dt 
\end{equation} 
Now, estimate the following inequality using Eqn. (37)
\begin{equation} \label{GrindEQ__3_15_} 
\frac{{}_{p} \Gamma _{k}^{\left(n-l\right)} \left(x\right)+_{p} \Gamma _{k}^{\left(n+l\right)} \left(x\right)}{2} -_{p} \Gamma _{k}^{\left(n\right)} \left(x\right) 
\end{equation} 
\begin{equation} \label{GrindEQ__3_16_} 
=\frac{1}{2} \left(\int _{0}^{\infty }t^{x-1} e^{-\frac{t^{k} }{p} } \log ^{n-l} t dt+\int _{0}^{\infty }t^{x-1} e^{-\frac{t^{k} }{p} } \log ^{n+l} t dt\right)-\int _{0}^{\infty }t^{x-1} e^{-\frac{t^{k} }{p} } \log ^{n} t dt 
\end{equation} 
\begin{equation} \label{GrindEQ__3_17_} 
=\frac{1}{2} \int _{0}^{\infty }\left(\frac{1}{\log ^{l} t} +\log ^{l} t-2\right)t^{x-1} e^{-\frac{t^{k} }{p} } \log ^{n} t dt\ge 0 
\end{equation} 

which implies
\begin{equation} \label{GrindEQ__3_18_} 
\frac{{}_{p} \Gamma _{k}^{\left(n-l\right)} \left(x\right)+_{p} \Gamma _{k}^{\left(n+l\right)} \left(x\right)}{2} \ge _{p} \Gamma _{k}^{\left(n\right)} \left(x\right) .                                                       
\end{equation} 

Exponentiating the above inequality yields the desired result.

  \textbf{Theorem 7.} \textit{For $x>0$and even integers $n\ge 1$, we have}

  \textit{i) }$e^{\psi _{k}^{\left(n\right)} \left(x\right)} \ge \sqrt{e^{\psi _{k}^{\left(n+1\right)} \left(x\right)} e^{\psi _{k}^{\left(n-1\right)} \left(x\right)} } $\textit{if n is odd and}

  \textit{ii) }$e^{\psi _{k}^{\left(n\right)} \left(x\right)} \ge \sqrt{e^{\psi _{k}^{\left(n+1\right)} \left(x\right)} e^{\psi _{k}^{\left(n-1\right)} \left(x\right)} } $\textit{if n is even.}

  \textit{Proof: }Estimate the following inequality using Eqn. \eqref{GrindEQ__1_5_}
\begin{equation} \label{GrindEQ__3_19_} 
\psi _{k}^{\left(n\right)} \left(x\right)-\frac{\psi _{k}^{\left(n+1\right)} \left(x\right)+\psi _{k}^{\left(n-1\right)} \left(x\right)}{2}  
\end{equation} 
\begin{equation} \label{GrindEQ__3_20_} 
=\left(-1\right)^{n+1} \left(\int _{0}^{\infty }\frac{t^{n} e^{-xt} }{1-e^{-kt} }  dt+\frac{1}{2} \int _{0}^{\infty }\frac{t^{n+1} e^{-xt} }{1-e^{-kt} }  dt+\frac{1}{2} \int _{0}^{\infty }\frac{t^{n-1} e^{-xt} }{1-e^{-kt} }  dt\right) 
\end{equation} 
\begin{equation} \label{GrindEQ__3_21_} 
=\frac{\left(-1\right)^{n+1} }{2} \int _{0}^{\infty }\frac{t^{n-1} e^{-xt} }{1-e^{-kt} }  \left(t+1\right)^{2} dt 
\end{equation} 

Thus, for odd n, we have
\begin{equation} \label{GrindEQ__3_22_} 
\psi _{k}^{\left(n\right)} \left(x\right)\ge \frac{\psi _{k}^{\left(n+1\right)} \left(x\right)+\psi _{k}^{\left(n-1\right)} \left(x\right)}{2}  
\end{equation} 

and for even n we have
\begin{equation} \label{GrindEQ__3_23_} 
\psi _{k}^{\left(n\right)} \left(x\right)\le \frac{\psi _{k}^{\left(n+1\right)} \left(x\right)+\psi _{k}^{\left(n-1\right)} \left(x\right)}{2}  
\end{equation} 
Exponentiating the above inequality yields the desired result.

\section{Conclusion}
In this paper, we derived some Turan-type inequalities for \textit{k-}polygamma function and \textit{p-k-}polygamma function using modified Holder's inequalities. The methodology used in this paper can also be applied to some other modified special functions such as the Nielsen's beta function and its generalizations and similar Turan-type inequalities can be obtained.

\end{document}